\documentclass[reqno,a4wide,english]{amsart}

\usepackage{amsmath}

\usepackage[all]{xy}
\CompileMatrices
\newdir{ >}{{}*!/-10pt/\dir{>}}

\theoremstyle{plain}
\newtheorem{theorem}{Theorem}[section]

\theoremstyle{definition}

\newtheorem{emptythm}[theorem]{}

\newcommand{\OXx}{{\mathcal O}_{X,x}}

\DeclareMathOperator{\Gl}{\mathrm{GL}}

\newcommand{\too}{\longrightarrow}
\newcommand{\sheaf}[1]{{\mathcal #1}}

\DeclareMathOperator{\Hom}{Hom}

\DeclareMathOperator{\Gal}{Gal}

\DeclareMathOperator{\Ker}{Ker}

\DeclareMathOperator{\Res}{Res}

\DeclareMathOperator{\id}{id}

\DeclareMathOperator{\khar}{char}
\DeclareMathOperator{\Spec}{Spec}

\DeclareMathOperator{\Norm}{N}

\DeclareMathOperator{\Kt}{\mathrm{K}}                       

\newcommand{\MK}{\Kt^{M}}                                          
\newcommand{\CM}{\mathrm{M}}                                  
\newcommand{\grAb}{\mathfrak{grAb}}                         
\newcommand{\MWK}{\mathrm{K}^{\mathrm{MW}}}    

\DeclareMathOperator{\W}{\mathrm{W}}

\DeclareMathOperator{\FdI}{\mathrm{I}}

\DeclareMathOperator{\Orth}{\mathrm{O}}

\DeclareMathOperator{\SymAlg}{\mathrm{S}}

\newcommand{\gk}{\mathrm{k}_{0}}        

\newcommand{\Fields}{\mathfrak{F}}
\newcommand{\bigFields}{\mathrm{Fields}}

\DeclareMathOperator{\Inv}{\mathrm{Inv}}  

\DeclareMathOperator{\HM}{\mathrm{H}}

\newcommand{\Z}{\mathbb{Z}}

\newcommand{\A}{\mathbb{A}}

\newcommand{\cf}{\textsl{cf.}\ }
\newcommand{\eg}{\textsl{e.g.}\ }
\newcommand{\ie}{\textsl{i.e.}\ }

\begin{document}

\title[On the splitting principle]
{On the splitting principle for cohomological invariants of reflection groups}

\author{Stefan Gille}
\email{gille@ualberta.ca}
\address{Department of Mathematical and Statistical Sciences,
University of Alberta, Edmonton T6G 2G1, Canada}
\author{Christian Hirsch}
\email{hirsch@uni-mannheim.de}
\address{Institut f\"ur Mathematik, Universit\"at Mannheim, B6 26,
68161 Mannheim, Germany}

\thanks{S.\ G.\ was supported by an NSERC grant.}
\thanks{C.\ H.\  was supported by The Danish Council for Independent Research |
Natural Sciences, grant DFF -- 7014-00074 \emph{Statistics for point processes in
space and beyond}, and by the \emph{Centre for Stochastic Geometry and
Advanced Bioimaging}, funded by grant 8721 from the Villum Foundation.}

\subjclass[2000]{Primary: 19D45; Secondary: 20G10}
\keywords{Reflection groups, Milnor-$K$-groups and cycle modules, cohomological invariants}

\date{September 17, 2019}

\begin{abstract}
Let~$\gk$ be a field and $W$ a finite orthogonal reflection group over~$\gk$.  We prove Serre's
splitting principle for cohomological invariants of~$W$ with values in Rost's cycle modules
(over~$\gk$) if the characteristic of~$\gk$ is coprime to~$|W|$. We then show that this principle
for such groups holds also for Witt- and Milnor-Witt $K$-theory invariants.
\end{abstract}

\maketitle

\section{Introduction}
\label{IntroSect}\bigbreak

\noindent
Let $\Fields_{\gk}$ be the category of finitely generated field extensions of a field~$\gk$,
and~$\CM_{\ast}$ a cycle module in the sense of Rost~\cite{Ro96} over~$\gk$, as for instance Milnor
$K$-theory or (abelian) Galois cohomology with finite coefficients. A {\it cohomological invariant} of
degree~$n$ of an algebraic group~$G$ over~$\gk$ with values in~$\CM_{\ast}$  is a natural
transformation
$$
a\, :\;\HM^{1}(\, -\, ,G)\,\too\,\CM_{n}(\, -\, )
$$
of functors on~$\Fields_{\gk}$. Here $\HM^{1}(\, -\, ,G)$ denotes the first non abelian Galois cohomology
of~$G$. Cohomological invariants are an old topic. For instance the discriminant, or the Clifford invariant
of a quadratic form can be interpreted as cohomological invariants of an orthogonal group. However the
formalization of this concept has been done only recently by Serre, see his lecture notes in~\cite{CohInv}
for a thorough account and some information on the history of the subject.

\smallbreak

In general the cohomological invariants of an algebraic group with values in a given cycle module are hard (if
not impossible) to compute. For most groups we know only some of the invariants, and even finding new ones
can be quit a task, as is exemplified in the construction of the Rost invariant, see \eg Merkurjev's lecture
in~\cite{CohInv}. Beside (the natural) applications to the classification of algebraic groups and their torsors there
are further applications of the theory of cohomological invariants, as for instance to rationality questions around
Noether's problem, see \eg~\cite[Part I, Sects.\ 33 and 34]{CohInv}.

\medbreak

The aim of this work, which is split into two parts, is the computation of the invariants of Weyl groups
with values in a cycle module, which is annihilated by~$2$, over a field of characteristic zero, which
contains a square root of~$-1$. The actual computation will be presented in the sequel~\cite{Hi19} to
this paper by the second named author, to which we also refer for a more precise description of the
result. Crucial for these investigations is the so called {\it splitting principle} for invariants of orthogonal
reflection groups. The proof of this principle is the content of this article. To formulate this result we recall
first the definition of a orthogonal reflection group over~$\gk$. Assume that $\khar\gk\not=2$. Let~$(V,b)$
be a regular symmetric bilinear space (of finite dimension) over~$\gk$ and~$\Orth (V,b)$ its orthogonal
group. A finite subgroup~$W$ of~$\Orth (V,b)$ is called a (finite) {\it orthogonal reflection group} over~$\gk$
if~$W$ is generated by reflections.

\medbreak

\noindent
{\bf Theorem.}
{\it
Let~$W$ be a orthogonal reflection group over the field~$\gk$. Assume
that $\khar {\gk}$ is coprime to the order of~$W$.
Then a cohomological invariant of degree~$n$ of~$W$ with values in a cycle
module~$\CM_{\ast}$ over~$\gk$
$$
a\, :\;\HM^{1}(\, -\, ,W)\,\too\,\CM_{n}(\, -\, )
$$
is trivial if and only if its restrictions to all $2$-subgroups of~$W$, which are generated
by reflections, are trivial.
}

\medbreak

\noindent
Crucial for our proof of this theorem is the explicit description of a versal $W$-torsor over~$\gk$,
which we give in Section~\ref{ReflGrVerTorSubSect}. This construction uses the fact that~$W$ is
a subgroup of the orthogonal group of some regular symmetric bilinear space over~$\gk$. Hence,
although~$W$ is a defined as a finite group scheme over an arbitrary field, we prove the
splitting principle only for invariants over fields, where~$W$ has a faithful orthogonal representation.

\smallbreak

We want to point out that our arguments here work also for Witt- and Milnor-Witt $K$-theory invariants of
orthogonal reflection groups, \ie the splitting principle holds also for such invariants. We explain the
necessary modifications in the last section of this work. However, the computations of Witt- and Milnor-Witt
$K$-theory invariants of Weyl groups are a different story and not touched in the second part of this
work. In fact, at least the Milnor-Witt $K$-theory invariants even of symmetric groups seem to be unknown.

\medbreak

For~$W$ a Weyl group the splitting principle in our theorem above has been already announced
by Serre in his lectures~\cite[Part I, 25.15]{CohInv}. It plays an important role in his computation of
cohomological invariants of the symmetric group with values in $\HM^{\ast}(\, -\, ,\Z/2)$ over an
arbitrary field. Note also that Ducoat claims this principle for the special case of invariants of Coxeter
groups with values in Galois cohomology with finite coefficients over (big enough) fields of characteristic
zero in his unpublished preprint~\cite{Du11}.

\medbreak

This article (except for the last section) as well as its sequel~\cite{Hi19} are based on the
2010 Diploma thesis~\cite{Hi10} of the second named author. This diploma thesis
does not deal with cycle modules of Rost, but with $\A^{1}$-invariants sheaves with
$\MK_{\ast}/2$-structure. However he proof of the splitting principle there has some
gaps and flaws.

\smallbreak

Our original intention was to write this article and its sequel in
the same setting but we refrained from this for the following two reasons. On the one hand,
it has turned out to be much easier and also shorter to give a complete proof of the splitting
principle in the slightly more restrictive setting of Rost's cycle modules. And on the other hand, we
believe that the most interesting invariants are anyway Galois cohomology-, or Milnor
$K$-theory (modulo some integer) invariants, which are both cycle modules, or Witt invariants.
An advantage of this restriction is also that the article is readable for readers only interested
in such invariants. They can assume throughout that the cycle module~$\CM_{\ast}$ in
question is one of their favorite theories.

\smallbreak

Moreover, according to Morel~\cite[Rem.\ 2.5]{A1AlgTop} such $\A^{1}$-invariant sheaves with
$\MK_{\ast}/2$-structure are the same as Rost cycle modules which are annihilated by~$2$, and
so we recover the computations of the diploma thesis~\cite{Hi10} in their full generality working
only with cycle modules.

\bigbreak

The content of this article is as follows. The main theorem is proven in Section~\ref{SpPrincipleSect}.
In Section~\ref{CycleModSect} we recall some definitions and facts about
Rost's cycle modules and Galois cohomology, mainly to fix notations and conventions. In the
following Section~\ref{WKSect} we recall the definition of a cohomological
invariant with values in a cycle module and prove some auxiliary results needed for
the proof of the main theorem. Except for the proof of a kind of specialization theorem,
our Theorem~\ref{specializationThm}, all arguments in Section~\ref{WKSect} are only
slight modifications of the one of Serre in his lectures~\cite[Part I]{CohInv}. However
the proof of Theorem~\ref{specializationThm} is more involved since a cycle module
over a field~$\gk$ is only defined for finitely generated field extensions of~$\gk$ and
so not for the completion or henselization of~$\gk$.

\smallbreak

Finally, in the last section we discuss the case of Witt-, and Milnor-Witt $K$-theory invariants.
Our proof of the splitting principle for invariants of reflection groups with values in cycle modules
carries over to this situation as well.

\bigbreak

\noindent
{\bf Acknowledgement.}
We would like to thank Fabien Morel for advise and fruitful discussions around this work.
This work has started (and slept long in between) more than 10 years ago, when one of us
(S.G.) was Assistent and the other (C.H.) Diploma Student of Fabien at the LMU Munich.
The first named author would also like to thank Volodya Chernousov, now his colleague
at the University of Alberta. Visiting Volodya in March 2008 has made a crucial impact on
this work.

\bigbreak\bigbreak

\begin{emptythm}
\label{NotationsSubSect}
{\bf Notations.}
Given a field~$\gk$ we denote by~$(\gk)_{s}$ its separable closure, and by
$\Gamma_{\gk}:=\Gal ((\gk)_{s}/\gk)$ its absolute Galois group.

\smallbreak

We denote by $\bigFields_{\gk}$ the category of all field
extensions of~$\gk$. More precisely, the objects of $\bigFields_{\gk}$ are pairs $(L,j)$,
where~$L$ is a field and $j:\gk\too L$ a homomorphism of fields. A morphism
$(E,i)\too (L,j)$ is a morphism of fields $\varphi:E\too L$,
such that
$$
\xymatrix{
E \ar[rr]^-{\varphi} & & L
\\
 & \gk \ar[ru]_-{j} \ar[lu]^-{i}
}
$$
commutes. For ease of notation the structure morphism will not be mentioned,
\ie we write~$L$ only instead of~$(L,j)$.

\smallbreak

If~$(\ell ,\iota)\in\bigFields_{\gk}$ then $\bigFields_{\ell}$ can be identified with a full subcategory
of~$\bigFields_{\gk}$ via the embedding $(L,j)\mapsto (L,j\circ\iota)$, which depends on
the structure morphism $\iota:\gk\too\ell$.

\smallbreak

The symbol~$\Fields_{\gk}$ denotes the full subcategory of~$\bigFields_{\gk}$ consisting of
finitely generated field extensions of~$\gk$, \ie of pairs~$(L,j)$, where~$L$ is a field and
$j:\gk\too L$ a homomorphism of fields giving~$L$ the structure of a finitely generated field
extension of~$\gk$. Again we can identify $\Fields_{\ell}$ with a full subcategory of~$\Fields_{\gk}$
for all~$\ell\in\Fields_{\gk}$.
\end{emptythm}

\goodbreak
\section{Preliminaries: Cycle modules, Galois cohomology and torsors}
\label{CycleModSect}\bigbreak

\begin{emptythm}
\label{CMSubSect}
{\bf Cycle modules.}
These have been invented by Rost~\cite{Ro96}
to facilitate Chow group computations. We refer to this article for details and more
information.

\smallbreak

The prototype of a cycle module is Milnor $K$-theory, which has been introduced by Milnor~\cite{Mi69/70},
and which we denote by
$$
\MK_{\ast}(F)\, :=\;\bigoplus\limits_{n\geq 0}\MK_{n}(F)
$$
for a field~$F$. Recall that this is a graded ring and as abelian group generated
by the pure symbols $\{ x_{1},\ldots ,x_{n}\}\in\MK_{n}(F)$, where $x_{1},\ldots ,x_{n}$
are non zero elements of~$F$.

\medbreak

A {\it cycle module over a field~$\gk$} is a covariant functor
$$
\CM_{\ast}\, :\;\Fields_{\gk}\,\too\,\grAb\, ,\; F\,\longmapsto\,\CM_{\ast}(F)\, =\,
\bigoplus\limits_{n\in\Z}\CM_{n}(F)\, ,
$$
where $\grAb$ denotes the category of graded abelian groups, such that $\CM_{\ast}(F)$
is a graded $\MK_{\ast}(F)$-module for all $F\in\Fields_{\gk}$. Following half way
Rost~\cite{Ro96} and deviating from somehow usual customs we denote by $\varphi_{\CM}$
the morphism $\CM_{\ast}(F)\too\CM_{\ast}(E)$ induced by a morphism $\varphi:F\too E$
in~$\Fields_{\gk}$.
\end{emptythm}

\begin{emptythm}
\label{2edResMapSubSect}
{\bf The second residue map.}
Let~$v$ be a discrete valuation of $F\in\Fields_{\gk}$ of geometric type which is trivial on~$\gk$.
By this we mean that there exists a normal integral $\gk$-scheme~$X$ of finite type, such that the function
field~$\gk (X)$ is equal~$F$, and such that~$v$ corresponds to a codimension one point in~$X$. Then
there is a $\MK_{\ast}(\gk)$-linear homomorphism, the so called {\it (second) residue map}:
$$
\partial_{v}\, :\;\CM_{\ast}(F)\,\too\,\CM_{\ast -1}(F(v))\, ,
$$
where~$F(v)$ is the residue field of~$v$.

\smallbreak

Associated with this homogenous homomorphism of degree~$-1$ there is a
homogenous homomorphism of degree~$0$, the so called {\it specialization
homomorphism}:
$$
s_{v}^{\pi}\, :\;\CM_{\ast}(F)\,\too\,\CM_{\ast}(F(v))\, ,\; x\,\longmapsto\,\partial_{v}\big(\{\pi\}\cdot x\big)\, ,
$$
which depends on the choice of a uniformizer~$\pi$ for~$v$.

\smallbreak

Recall the following three axioms, which play some role in the next section.
Let $F,v,F(v)$ be as above and $\varphi: F\too E$ a finite field extension. Assume
there is a geometric valuation~$w$ on~$E$ with residue field~$E(w)$ and with
$w|_{F}=v$. Let~$e_{w|v}$ be the ramification index and $\bar{\varphi}:F(v)\too E(w)$
the induced homomorphism of the residue fields. Then
the following holds (numbering as in Rost~\cite[p.\ 329]{Ro96}):

\smallbreak

\begin{itemize}
\item[{\bf (R3a)}]
$\partial_{w}\circ\varphi_{\CM}\, =\, e_{w|v}\cdot\bar{\varphi}_{\CM}\circ\partial_{v}$;

\smallbreak

\item[{\bf (R3c)}]
if~$w$ is trivial on~$F$, and so~$F(v)=F$, then $\partial_{w}\circ\varphi_{\CM}\, =0$; and

\smallbreak

\item[{\bf (R3d)}]
if~$w$ is as in {\bf (R3c)} and~$\pi$ is an uniformizer for~$w$ then $s_{w}^{\pi}\circ\varphi_{\CM}\, =\,\bar{\varphi}_{\CM}$.
\end{itemize}
\end{emptythm}

\begin{emptythm}
\label{UnrCMSubSect}
{\bf Unramified cycle modules.}
Let~$X$ be a integral scheme, which is essentially of finite type over a field~$\gk$. By
the latter we mean that~$X$ is a finite type $\gk$-scheme or a localization of such a
scheme. Denoting~$X^{(1)}$ the points of codimension~$1$. If a point~$x$ in~$X^{(1)}$
is regular, then its local ring~$\OXx$ is a discrete valuation ring and we get a valuation~$v_{x}$
on the function field~$\gk (X)$ of~$X$.

\smallbreak

Given a cycle module~$\CM_{\ast}$ we have then a second residue map
$$
\partial_{x}\, :\; =\,\partial_{v_{x}}\, :\;\CM_{\ast}(\gk (X))\,\too\,\CM_{\ast -1}(\gk (x))\, ,
$$
where~$\gk(x)$ denotes the residue field of~$x$, as well as a specialization map
$$
s_{x}^{\pi}\, :=\; s_{v_{x}}^{\pi}\, :\;\CM_{\ast}(\gk (X))\too\CM_{\ast}(\gk (x))
$$
for every uniformizer $\pi\in\OXx$.
If~$X$ is regular in codimension one $\partial_{x}$ exists for all~$x\in X^{(1)}$ and so
we can define
$$
\CM_{n,unr}(X)\, :=\;\Ker\big(\,\CM_{n}(\gk (X))\,\xrightarrow{\; (\partial_{x})_{x\in X^{(1)}}\;}\,
\bigoplus\limits_{x\in X^{(1)}}\CM_{n-1}(\gk (x))\, ,
$$
the so called {\it unramified $\CM_{n}$-cohomology group} of~$X$.
\end{emptythm}

\begin{emptythm}
\label{H1SubSect}
{\bf Non abelian Galois cohomology.}
We recall now -- mainly to fix notations -- some definitions and properties of torsors and
non abelian Galois cohomology sets. We refer to Serre's well known book~\cite{GalCoh}
and also to~\cite[\S 28 and \S 29]{BookInv} for details and more information.

\medbreak

Let~$F$ be a field and $G$ a linear algebraic group over~$F$. We denote by
$\HM^{1}(F,G)$ the {\it first non abelian Galois cohomology} set,
\ie $\HM^{1}(F,G)=\HM^{1}(\Gamma_{F},G(F_{s}))$. If a continuous maps
$c:\Gamma_{F}\too G(F_{s})$, $\sigma\mapsto c_{\sigma}$ is a cycle, \ie represents
an element of~$\HM^{1}(F,G)$, then we denote its class in~$\HM^{1}(F,G)$ by~$[c]$.

\smallbreak

If $\varphi:F\too E$ is a morphism of fields we denote the induced {\it restriction map}
$\HM^{1}(F,G)\too\HM^{1}(E,G)$ by $r_{\varphi}$, or if~$\varphi$ is clear from the
context by $r_{E/F}$.

\smallbreak

If $\theta:H\too G$ is a morphism of linear algebraic groups over~$F$ we denote
following~\cite{BookInv} the induced homomorphism $\HM^{1}(F,H)\too\HM^{1}(F,G)$
by~$\theta^{1}$.

\smallbreak

In the proof of Theorem~\ref{specializationThm} below we consider also the first non
abelian \'etale cohomology set $\HM^{1}_{et}(X,G)$, where~$X$ is a scheme over~$F$
and~$G$ a linear algebraic group over~$F$. If $f:X\too Y$ is a morphism of such schemes
we denote the pull-back map $\HM^{1}_{et}(Y,G)\too\HM^{1}_{et}(X,G)$ by~$r_{f}$. We write
then also $T_{X}$ instead of $r_{f}(T)$ for $T\in\HM^{1}_{et}(Y,G)$ if~$f$ is clear from
the context.

\smallbreak

Note that since~$G$ is smooth the set~$\HM^{1}_{et}(X,G)$ can be identified with the isomorphism
classes of $G$-torsors $\pi:\sheaf{T}\too X$ over~$X$, see \eg~\cite[Chap.\ III.4]{EC}. We denote
the class of a $G$-torsor $\pi:\sheaf{T}\too X$ over~$X$ by $[\sheaf{T}\too X]$.

\smallbreak

We use in the following also affine notations, \ie $\HM^{1}_{et}(R,G)$ instead of $\HM^{1}_{et}(X,G)$
if $X=\Spec R$ is affine.  Note that if~$X=\Spec K$ is the spectrum of a field then
$\HM^{1}_{et}(X,G)$ is naturally isomorphic to~$\HM^{1}(K,G)$.
\end{emptythm}

\begin{emptythm}
\label{GaloisExpl}
{\bf Example.}
Let~$G$ be a finite group with trivial $\Gamma_{F}$-action, where~$F$ is a field.
Then the non abelian Galois cohomology set $\HM^{1}(F,G)$ can be identified with
the isomorphism classes of $G$-Galois algebras, see \eg~\cite[V.14]{IGC},
or~\cite[\S 18B]{BookInv}.

\smallbreak

A particular example of such a $G$-Galois algebra is a finite Galois extension $E\supset F$
with group $\Gal (E/F)=G$. Then the continuous homomorphism of groups
$c:\Gamma_{F}\too G$, $\sigma\mapsto\sigma|_{E}$, represents the class of the
$G$-Galois algebra~$E$ in~$\HM^{1}(F,G)$, see \eg~\cite[Thm.\ V.14.17]{IGC}.

\smallbreak

If $\theta:H\subset G$ is a subgroup with fixed field~$L$ then the
class of the restriction $r_{L/F}([c])$ is represented by the continuous homomorphism
$c|_{\Gamma_{L}}:\Gamma_{L}\subset\Gamma_{F}\xrightarrow{c}G$, whose
image is in the subgroup~$H$. It follows that
$$
r_{L/F}([c])\, =\,\theta^{1}([c'])\, ,
$$
where $[c']\in\HM^{1}(L,H)$ is represented by
$c':\Gamma_{L}\xrightarrow{c|_{\Gamma_{L}}}H$.
\end{emptythm}

\begin{emptythm}
\label{VerTorSubSect}
{\bf Versal torsors.}
Let~$T\in\HM^{1}_{et}(X,G)$ be a $G$-torsor over the smooth integral $F$-scheme~$X$
with function field $K=F(X)$. Assume that given a field extension~$L$ of~$F$, which is infinite,
and an element $\tilde{T}\in\HM^{1}(L,G)$ then there exists an $L$-point~$x$ of~$X$, such that
$T_{F(x)}=\tilde{T}$ in $\HM^{1}(F(x),G)=\HM^{1}(L,G)$. Then~$T_{K}\in\HM^{1}(K,G)$ is called
a {\it versal $G$-torsor}, see~\cite[Part I, Def.\ 5.1]{CohInv}.
\end{emptythm}

\begin{emptythm}
\label{VerTorExpl}
{\bf Example.}
Let~$G$ be a finite group which acts faithfully on the finite dimensional
$\gk$-vector space~$V$, where~$\gk$ is a field. Then~$G$ acts on the
dual space $V^{\vee}:=\Hom_{\gk}(V,\gk)$ via $(h.f)(v):=f(h^{-1}.v)$
for all $f\in V^{\vee}$, $v\in V$, and~$h\in G$. This induces a
$G$-action on $\A (V):=\Spec\SymAlg (V^{\vee})$, where~$\SymAlg (V^{\vee})$
denotes the symmetric algebra of~$V^{\vee}$. For~$g\in G$ denote by
$\sheaf{V}_{g}$ the closed subset of~$\A (V)$ defined by the ideal generated
by all $f\circ (g-\id_{V})\in V^{\vee}$, $f\in V^{\vee}$. The group~$G$ acts freely
on the open set
$$
U\, :=\;\A (V)\,\setminus\;\bigcup\limits_{g\in G}\sheaf{V}_{g}\, ,
$$
and so the quotient morphism $q:U\too U/G$ is a $G$-torsor. The generic
fiber of this torsor is a versal $G$-torsor, see~\cite[Part I, 5.4 and 5.5]{CohInv}
for a proof. Note that the function field~$\gk (U/G)$ is the fraction field of the
invariant ring $\SymAlg (V^{\vee})^{G}$. Hence the class of this versal $G$-torsor
in $\HM^{1}(\gk (U/G),G)$ is the class of the $G$-Galois algebra $\gk (U)$, \ie
the class of the Galois extension $\gk (U)\supseteq\gk(U)^{G}=\gk (U/G)$.
\end{emptythm}

\goodbreak
\section{Invariants in cycle modules}
\label{WKSect}\bigbreak

\begin{emptythm}
\label{InvCMDef}
{\bf Definition.}
Let~$G$ be a linear algebraic group and~$\CM_{\ast}$ a cycle module over the field~$\gk$.
A {\it cohomological invariant of degree~$n$} of~$G$ with values in the cycle module~$\CM_{n}$ is
a natural transformation of functors
$$
a\, :\;\HM^{1}(\, -\, ,G)\,\too\,\CM_{n}(\, -\,)\, ,
$$
\ie for all $\varphi:F\too E$ in~$\Fields_{\gk}$ the following diagram commutes:
$$
\xymatrix{
\HM^{1}(E,G) \ar[r]^-{a_{E}} & \CM_{n}(E)
\\
\HM^{1}(F,G) \ar[r]_-{a_{F}} \ar[u]^-{r_{\varphi}} & \CM_{n}(F) \ar[u]_-{\varphi_{\CM}} \rlap{\, .}
}
$$
This definition is due to Serre and a special case of the one given in his lectures~\cite[Part I]{CohInv}, but
it includes the main player of Serre's text, cohomological invariants of algebraic groups with coefficients
in Galois cohomology. Note however that there is a subtle difference as we consider here only the
category~$\Fields_{\gk}$ of finitely generated field extensions of~$\gk$ and not the category~$\bigFields_{\gk}$
of all field extensions of~$\gk$. This is forced by the fact that for technical reasons a "abstract" cycle module
is not defined for all field extensions of a given field but only for the finitely generated ones. If one is
only interested in "concrete" cycle modules as for instance Milnor $K$-theory, Witt groups, or Galois
cohomology, there is no need for this restriction. However one may wonder whether it could happen
albeit~$\HM^{1}(\, -\, ,G)$ and the value group is defined for all $F\in\bigFields_{\gk}$ that there exists
invariants, which are only defined for finitely generated field extensions of~$\gk$.

\medbreak

Following Serre's lectures~\cite[Part I]{CohInv} we denote the set of cohomological invariants of degree~$n$
of the group~$G$ with values in~$\CM_{n}$ by $\Inv_{\gk}^{n}(G,\CM_{\ast})$. The set $\Inv^{n}_{\gk}(G,\CM_{\ast})$
has the structure of an is an abelian group as~$\CM_{n}(F)$ is one for all $F\in\Fields_{\gk}$.

\smallbreak

We set
$$
\Inv_{\gk}(G,\CM_{\ast})\, :=\;
\bigoplus\limits_{n\in\Z}\Inv^{n}_{\gk}(G,\CM_{\ast})\, ,
$$
and call elements of this direct sum {\it cohomological invariants} of~$G$ with
values in~$\CM_{\ast}$. Note that the $\MK_{\ast}$-structure of~$\CM_{\ast}$
induces a $\MK_{\ast}(\gk)$ operation on~$\Inv_{\gk}(G,\CM_{\ast})$ making
the set of cohomological invariants of~$G$ with values in~$\CM_{\ast}$ a
graded $\MK_{\ast}(\gk)$-module.
\end{emptythm}

\begin{emptythm}
\label{ConstInvSubSect}
{\bf Constant invariants.}
We have $\Inv_{\gk} (G,\CM_{\ast})\not= 0$ if and only $\CM_{\ast}(\gk)\not= 0$. In fact,
if $x\in\CM_{\ast}(\gk)$ then
$$
c_{L}\, :\;\HM^{1}(L,G)\,\too\,\CM_{\ast}(L)\, ,\; t\,\longmapsto\, (\iota_{L})_{\CM}(x)\, ,
$$
where $L\in\Fields_{\gk}$ with structure map $\iota_{L}:\gk\too L$, defines an
invariant. Such invariants are called {\it constant}, respectively, {\it constant of degree~$n$}
if $x\in\CM_{n}(\gk)$, and we write $c\equiv x\in\CM_{\ast}(\gk)$.
\end{emptythm}

\begin{emptythm}
\label{ResInvSubSect}
{\bf Restriction of invariants.}
Let $\theta:H\too G$ be a morphism of linear algebraic groups over~$\gk$ and~$\CM_{\ast}$
a cycle module over~$\gk$. Composing $a\in\Inv_{\gk} (G,\CM_{\ast})$ with the map~$\theta^{1}$:
$$
\HM^{1}(\, -\, ,H)\,\xrightarrow{\;\theta^{1}\;}\HM^{1}(\, -\, ,G)\,\xrightarrow{\; a\;}\,\CM_{\ast}(\, -\,)
$$
we get an invariant $\theta^{\ast}(a)\in\Inv_{\gk} (H,\CM_{\ast})$. In case $\theta:H\subseteq G$
is the embedding of a closed subgroup we denote following Serre's lecture~\cite[Part I]{CohInv}
the induced homomorphism $\theta^{\ast}:\Inv_{\gk} (G,\CM_{\ast})\too\Inv_{\gk} (H,\CM_{\ast})$ by
$\Res_{G}^{H}$ and call it the {\it restriction}.
\end{emptythm}

\begin{emptythm}
\label{InnAutExpl}
{\bf Example.}
Let~$H$ be a subgroup of a finite group~$G$, $\CM_{\ast}$ a cycle
module over~$\gk$, and $a\in\Inv_{\gk} (G,\CM_{\ast})$. If~$g$ is an
element of the normalizer~$\Norm_{G}(H)$ we denote by~$\iota_{g}$
the inner automorphism of~$G$ defined by~$g$. The isomorphism $\iota_{g}$ acts
also on~$H$ and so consequently via $\iota_{g}^{\ast}:a\mapsto a\circ\iota_{g}^{1}$
on $\Inv_{\gk} (H,\CM_{\ast})$ for all $g\in\Norm_{G}(H)$ giving~$\Inv_{\gk} (H,\CM_{\ast})$ the
structure of a $\Norm_{G}(H)$-module.  We claim that $\Res_{G}^{H}$ maps
$\Inv_{\gk} (G,\CM_{\ast})$ into the subgroup $\Inv_{\gk} (H,\CM_{\ast})^{\Norm_{G}(H)}$ of
$\Norm_{G}(H)$-invariant elements in~$\Inv_{\gk} (H,\CM_{\ast})$. In fact, by~\cite[Part I, Prop.\ 13.1]{CohInv}
the map $\iota_{g}^{1}:\HM^{1}(L,G)\too\HM^{1}(L,G)$ is the identity for all $L\in\Fields_{\gk}$.
Since we have $\Res_{G}^{H}\circ\iota_{g}^{\ast}=\iota_{g}^{\ast}\circ\Res_{G}^{H}$ this implies the claim.

\medbreak

We prove now the following specialization theorem, which is kind of an analog of a result
of Rost~\cite[Part I, Thm.\ 11.1]{CohInv} about Galois cohomology invariants.
\end{emptythm}

\begin{emptythm}
\label{specializationThm}
{\bf Theorem.}
{\it
Let~$X$ be an integral scheme with function field~$K=\gk(X)$, which is essentially of finite type
over the field~$\gk$. Let further~$\CM_{\ast}$ be a cycle module over~$\gk$, $a\in\Inv_{\gk} (G,\CM_{\ast})$,
where~$G$ is a linear algebraic group over~$\gk$, and $T\in\HM^{1}_{et}(X,G)$. Let~$x\in X$
be a regular codimension one point. Then we have:

\smallbreak

\begin{itemize}
\item[(i)]
$\partial_{x}\big(\, a_{K}(T_{K})\,\big)\, =0$; and

\smallbreak

\item[(ii)]
$s_{x}^{\pi}(a_{K}(T_{K}))\, =\, a_{\gk (x)}(T_{\gk (x)})$ for all local uniformizers~$\pi\in\OXx$.
\end{itemize}
In particular, if~$X$ is regular in codimension one, then $a_{K}(T_{K})\in\CM_{\ast,unr}(X)$.
}

\begin{proof}
Replacing~$X$ by~$\Spec\OXx$ we can assume that $X=\Spec R$ for a discrete valuation
ring~$R$, which is essentially of finite type over~$\gk$, and that~$x$ is the closed point of~$X$.
We denote~$k=\gk (x)$ the residue field of~$R$, $q:R\too k$ the quotient map, $\eta:R\too K$ the embedding
of~$R$ into its field of fractions~$K$, and~$v$ the discrete valuation of~$K$ corresponding to~$R$.
Then $\partial_{x}=\partial_{v}:\CM_{\ast}(K)\too\CM_{\ast -1}(k)$, and for~(i) we have to show
$\partial_{v}\big(\, a_{K}(T_{K})\,\big)=0$.

\smallbreak

To prove this let $\psi: R\too R^{h}$ be the henselization of~$R$ with fraction field~$K^{h}$.
This is also a discrete valuation ring with the same residue field~$k$, and there exists local
\'etale $R$-algebras $\varphi_{i}:R\too R_{i}$, $i\in I$, such that $R^{h}=\lim\limits_{i\in I} R_{i}$,
see~\cite[Chap.\ VIII]{ALH}. The rings~$R_{i}$ are also discrete valuation rings with~$k$ as
residue field. We denote by~$v_{i}$ the induced valuation on the fraction field~$K_{i}$
of~$R_{i}$. We get a commutative diagram of homomorphisms of rings for all~$i\in I$:
$$
\xymatrix{
K \ar[r]^-{\varphi'_{i}} & K_{i} \ar[r]^-{\psi'_{i}} & K^{h}
\\
R \ar[r]_-{\varphi_{i}} \ar[u]^-{\eta} & R_{i} \ar[r]_-{\psi_{i}} \ar[u]_-{\eta_{i}} & R^{h} \ar[u]_-{\eta^{h}} \rlap{\, ,}
}
$$
where the up going arrows are the respective inclusions of the rings~$R,R_{i}$, and~$R^{h}$ into
their fraction fields.

\medbreak

The homomorphism $\varphi_{i}:R\too R_{i}$ is unramified at the maximal ideal of~$R_{i}$, and so
by the cycle module Axiom~{\bf (R3a)}, see Section~\ref{2edResMapSubSect}, the square on the
right hand side of the following diagram commutes
\begin{equation}
\label{SpPfEq0}
\xymatrix{
\HM^{1}(K_{i},G) \ar[r]^-{a_{K_{i}}} & \CM_{\ast}(K_{i}) \ar[r]^-{\partial_{v_{i}}} & \CM_{\ast -1}(k)
\\
\HM^{1}(K,G) \ar[r]^-{a_{K}} \ar[u]^-{r_{\varphi_{i}'}} & \CM_{\ast}(K) \ar[u]^-{(\varphi_{i}')_{\CM}}
     \ar[r]^-{\partial_{v}} & \CM_{\ast -1}(k) \ar[u]_-{=}
}
\end{equation}
for all~$i\in I$. Since~$a$ is an invariant also the one on the left hand side is commutative.
Therefore to prove $\partial_{v}(a_{K}(T_{K}))=0$ it is enough to show that there exists~$i\in I$,
such that
$$
\partial_{v_{i}}\big(\, a_{K_{i}}(T_{K_{i}})\,\big)\, =\,
\partial_{v_{i}}\big(\,a_{K_{i}}(r_{\varphi'_{i}}(T_{K}))\,\big)\; =\, 0\, .
$$
To find this~$i\in I$ we use the fact that there exists a splitting $j: k\too R^{h}$
of the quotient morphism $q^{h}:R^{h}\too k$, \ie $q^{h}\circ j=\id_{k}$. In fact,
if~$\hat{R}$ is the completion of~$R$ we have a splitting $\hat{j}:k\too\hat{R}$
of the quotient map $\hat{R}\too k$ by the Cohen structure theorem. This splitting
factors via~$R^{h}$ since the henselian local domain~$R^{h}$ is excellent
by~\cite[Cor.\ 18.7.6]{EGA4-4}, and therefore satisfies the approximation property
by~\cite[Sect.\ 3.6, Cor.\ 9]{NM}, which implies in particular, that the splitting~$\hat{j}$
of $\hat{R}\too k$ factors via~$R^{h}$.

\smallbreak

The composition of maps
$\HM^{1}(k,G)\xrightarrow{r_{j}}\HM^{1}_{et}(R^{h},G)\xrightarrow{r_{q^{h}}}\HM^{1}(k,G)$
is the identity, and by~\cite[Chap.\ XXIV, Prop.\ 8.1]{SGA3-3} the map
$r_{q^{h}}:\HM^{1}_{et}(R^{h},G)\too\HM^{1}(k,G)$ is an isomorphism. Therefore
$r_{j}:\HM^{1}(k,G)\too\HM^{1}_{et}(R^{h},G)$ is one as well, and moreover we have
$r_{j}(T_{k})=T_{R^{h}}$ since $r_{q^{h}}(T_{R^{h}})=r_{q^{h}}(r_{\psi}(T_{R}))=r_{q}(T_{R})=T_{k}$.

\smallbreak

Let~$k_{i}$ be the pre-image of $j(k)$ under the homomorphism $\psi_{i}:R_{i}\too R^{h}$.
The set~$k_{i}\setminus\{ 0\}$ is contained in the set of units of~$R_{i}$ and so
is a field. This implies also that $v_{i}$ is trivial on~$k_{i}$. The $\gk$-linear quotient
homomorphism $q_{i}:R_{i}\too k$ maps~$k_{i}$ onto a subfield of~$k$. This implies
by~\cite[Chap.\ 5, \S 14, No 7, Cor.\ 3]{ALG4-7} that~$k_{i}$ is also a finitely generated
field extension of~$\gk$ and so in~$\Fields_{\gk}$ for all~$i\in I$.

\smallbreak

By the definition of the fields~$k_{i}$ we have a commutative diagram
\begin{equation}
\label{SpPfEq2}
\xymatrix{
k \ar[r]^-{j} & R^{h}
\\
k_{i} \ar[r]_-{j_{i}} \ar[u]^-{\bar{\psi}_{i}} & R_{i} \ar[u]_-{\psi_{i}}
}
\end{equation}
for all~$i\in I$, where $j_{i}$ is the inclusion $k_{i}\subset R_{i}$
and $\bar{\psi}_{i}$ the homomorphism induced by~$\psi_{i}$. Note that
$\bar{\psi}_{i}=q_{i}\circ j_{i}$ as $q\circ j=\id_{k}$.

\smallbreak

Diagram~(\ref{SpPfEq2}) gives in turn a commutative
diagram of pointed non abelian cohomology sets
\begin{equation}
\label{SpPfEq3}
\xymatrix{
\HM^{1}(k,G) \ar[r]^-{r_{j}} & \HM^{1}_{et}(R^{h},G)
\\
\HM^{1}(k_{i},G) \ar[r]_-{r_{j_{i}}} \ar[u]^-{r_{\bar{\psi}_{i}}} &
    \HM^{1}_{et}(R_{i},G) \ar[u]_-{r_{\psi_{i}}} \rlap{\, .}
}
\end{equation}

\smallbreak

We have $k=\lim\limits_{i\in I}k_{i}$ and therefore by~\cite[Chap.\ VII, Thm.\ 5.7]{SGA4-2}
(or by direct verification) that $\HM^{1}(k,G)\, =\,\lim\limits_{i\in I}\HM^{1}(k_{i},G)$. Hence there
exists $i_{0}\in I$ and $T_{i_{0}}\in\HM^{1}(k_{i_{0}},G)$, such that
\begin{equation}
\label{SpPfEq4}
r_{\bar{\psi}_{i_{0}}}(T_{i_{0}})\, =\, T_{k}\in\HM^{1}(k,G)\, .
\end{equation}

\medbreak

By~(\ref{SpPfEq3}) we have
$$
r_{\psi_{i_{0}}}\big(\, r_{j_{i_{0}}}(T_{i_{0}})\,\big)\, =\,
r_{j}\big(\, r_{\bar{\psi}_{i_{0}}}(T_{i_{0}})\,\big)\, =\, r_{j}(T_{k})\, =\, T_{R^{h}}\, =\,
r_{\psi_{i_{0}}}(T_{R_{i_{0}}})\, .
$$
Now $R^{h}=\lim\limits_{i\in I}R_{i}$ and so by~\cite[Chap.\ VII, Thm.\ 5.7]{SGA4-2} again we
have $\lim\limits_{i\in I}\HM^{1}_{et}(R_{i},G)=\HM^{1}_{et}(R^{h},G)$. Hence replacing~$i_{0}$
by a 'larger' element of~$I$ if necessary we can assume that also
\begin{equation}
\label{SpPfEq5}
r_{j_{i_{0}}}(T_{i_{0}})\, =\, T_{R_{i_{0}}}\, .
\end{equation}

\medbreak

We claim that this index~$i_{0}$ does the job, \ie we have
$\partial_{v_{i_{0}}}\big( a_{K_{i_{0}}}(T_{K_{i_{0}}})\big)=0$.
In fact, since~$a$ is a cohomological invariant we
have a commutative diagram
$$
\xymatrix{
\HM^{1}(K_{i_{0}},G) \ar[r]^-{a_{K_{i_{0}}}} & \CM_{\ast}(K_{i_{0}})
\\
\HM^{1}(k_{i_{0}},G) \ar[r]_-{a_{k_{i_{0}}}} \ar[u]^-{r_{(\eta_{i_{0}}\circ j_{i_{0}})}} &
    \CM_{\ast}(k_{i_{0}}) \ar[u]_-{(\eta_{i_{0}}\circ j_{i_{0}})_{\CM}} \rlap{\, ,}
}
$$
and therefore taking~(\ref{SpPfEq5}) into account
$$
a_{K_{i_{0}}}(T_{K_{i_{0}}})\, =\, a_{K_{i_{0}}}\big(\, r_{(\eta_{i_{0}}\circ j_{i_{0}})}(T_{i_{0}})\,\big)\, =\,
(\eta_{i_{0}}\circ j_{i_{0}})_{\CM}\big(\, a_{k_{i_{0}}}(T_{i_{0}})\,\big)\, .
$$
But $v_{i_{0}}|_{k_{i_{0}}}\equiv 0$ and so by the cycle module Axiom~{\bf (R3c)}, see
Section~\ref{2edResMapSubSect}, we have $\partial_{v_{i_{0}}}(z)=0$ for all $z\in\CM_{\ast}(K_{i_{0}})$,
which are in the image of $(\eta_{i_{0}}\circ j_{i_{0}})_{\CM}$. We have proven~(i).

\bigbreak

For the proof of~(ii) we continue with above notation, \ie $R=\OXx$,
$R^{h}=\lim\limits_{i\in I}R_{i}$, and so on. We fix further a uniformizer~$\pi$ of~$R$.
Since the extensions $\varphi_{i}:R\too R_{i}$ are unramified the image of the uniformizer~$\pi$
in~$R_{i}$ is also one, which we denote also by~$\pi$. We have $s_{x}^{\pi}=s_{v}^{\pi}$,
and so taking the right hand side of the commutative diagram~(\ref{SpPfEq0}) as well as the
definition of the specialization map, see Section~\ref{2edResMapSubSect}, into account we have
$s_{v_{i_{0}}}^{\pi}\circ (\varphi'_{i_{0}})_{\CM}\, =\, s_{v}^{\pi}$. Hence we have
$$
\begin{array}{r@{\; =\;}l@{\qquad}l}
      s_{v}^{\pi}(a_{K}(T_{K})) & s_{v_{i_{0}}}^{\pi}\big(\, (\varphi'_{i_{0}})_{\CM}(a_{K}(T_{K}))\,\big) & \\[4mm]
           & s_{v_{i_{0}}}^{\pi}\big(\, a_{K_{i_{0}}}(r_{\varphi'_{i_{0}}}(T_{K}))\,\big) & \mbox{$a$ is invariant} \\[4mm]
           & s_{v_{i_{0}}}^{\pi}\big(\, a_{K_{i_{0}}}(r_{\eta_{i_{0}}}(T_{R_{i_{0}}}))\,\big) &
                  \mbox{since $\eta_{i_{0}}\circ\varphi_{i_{0}}=\varphi'_{i_{0}}\circ\eta$} \\[4mm]
           & s_{v_{i_{0}}}^{\pi}\big(\, a_{K_{i_{0}}}(r_{(\eta_{i_{0}}\circ j_{i_{0}})}(T_{i_{0}}))\big) &
                  \mbox{by~(\ref{SpPfEq5})} \\[4mm]
           & s_{v_{i_{0}}}^{\pi}\big(\, (\eta_{i_{0}}\circ j_{i_{0}})_{\CM}(a_{k_{i_{0}}}(T_{i_{0}}))\,\big) &
                 \mbox{$a$ is invariant} \\[4mm]
           & (\bar{\psi}_{i_{0}})_{\CM}\big(\, a_{k_{i_{0}}}(T_{i_{0}})\,\big) & \mbox{by~{\bf (R3d)}} \\[4mm]
           & a_{k}\big(\, r_{\bar{\psi}_{i_{0}}}(T_{i_{0}})\,\big) & \mbox{$a$ is invariant} \\[4mm]
           & a_{k}(T_{k}) & \mbox{by (\ref{SpPfEq4}).}
\end{array}
$$
as claimed. We are done.
\end{proof}

\smallbreak

A consequence of this result is the following detection principle, which
is the analog of~\cite[Part I,12.2]{CohInv} for cycle module invariants.
\end{emptythm}

\begin{emptythm}
\label{SpCor}
{\bf Corollary.}
{\it
Let~$R$ be a regular local ring, which is essentially of finite type over the field~$\gk$.
Denote by~$K$ and~$k$ the fraction- and residue field, respectively, of~$R$. Let
further~$G$ be a linear algebraic group over~$\gk$, and~$\CM_{\ast}$ a cycle module
over~$\gk$. Then
$$
a_{K}(T_{K})\, =\, 0\quad\Longrightarrow\quad a_{k}(T_{k})\, =\, 0\
$$
for all $T\in\HM^{1}_{et}(R,G)$ and all $a\in\Inv_{\gk}(G,\CM_{\ast})$.
}

\begin{proof}
The proof is the same as the one of~\cite[Part I, 12.2]{CohInv}. We recall
the arguments for the convenience of the reader.

\smallbreak

If $\dim R=1$ this follows from part~(ii) of the theorem above, so let~$d:=\dim R\geq 2$,
and~$t\in R$ a regular parameter. Then $R/Rt$ is also a regular local ring with the same
residue field~$k$, and which is essentially of finite type over~$\gk$. The quotient field
of~$R/Rt$ is the residue field~$K_{t}$ of the discrete valuation ring $R_{Rt}$ (the localization at the
codimension one prime ideal~$Rt$). By the dimension one case we have $a_{K_{t}}(T_{K_{t}})=0$,
and so by induction $a_{k}(T_{k})=0$.
\end{proof}

\medbreak

Finally we state and prove the following detection principle, which is the cycle module
analog of~\cite[Part I, Thm.\ 12.3]{CohInv}. Again the proof is the same as in Serre's lecture
and only recalled for the convenience of our reader.
\end{emptythm}

\begin{emptythm}
\label{DetectionThm}
{\bf Theorem.}
{\it
Let~$G$ be a linear algebraic group over the field~$\gk$, and $T\in\HM^{1}(K,G)$
a versal $G$-torsor. Then we have for a given cycle module~$\CM_{\ast}$ over~$\gk$
and $a,b\in\Inv_{\gk} (G,\CM_{\ast})$:
$$
a_{K}(T)\, =\, b_{K}(T)\quad\Longrightarrow\quad a\, =\, b\, .
$$
}

\begin{proof}
Replacing~$a$ by $b-a$ it is enough to show that $a_{K}(T)=0$ implies~$a\equiv 0$.
We have to show $a_{k}(S)=0$ for all $k\in\Fields_{\gk}$ and all $S\in\HM^{1}(k,G)$.

\smallbreak

Replacing~$k$ by the rational function field~$k(T)$ if necessary, we can assume
that~$k$ is an infinite field. In fact, since~$a$ is an invariant the following diagram
commutes:
$$
\xymatrix{
\HM^{1}(k(T),G) \ar[r]^-{a_{k(T)}} & \CM_{\ast}(k(T))
\\
\HM^{1}(k,G) \ar[r]^-{a_{k}} \ar[u]^-{r_{\iota}} & \CM_{\ast}(k)
    \ar[u]_-{\iota_{\CM}} \rlap{\, ,}
}
$$
where $\iota:k\hookrightarrow k(T)$ is the natural embedding. By
Rost~\cite[Prop.\ 2.2 {\bf (H)}]{Ro96} the homomorphism
$\iota_{\CM}:\CM_{\ast}(k)\too\CM_{\ast}(k(T))$ is injective and so $a_{k(T)}(S_{k(T)})=0$
implies $a_{k}(S)=0$.

\smallbreak

To prove the claim for an infinite field~$k$ we use that since~$T\in\HM^{1}(K,G)$ is a versal
$G$-torsor there exists a $G$-torsor $\sheaf{T}\too X$ over a smooth integral scheme~$X$ with
function field~$K$, such that the generic fiber of $\sheaf{T}\too X$ is isomorphic to~$T$,
and such that there exists $x\in X(k)$ with $S=[\sheaf{T}\too X]_{\gk (x)}$. Now~$\OXx$ is a
regular local ring since~$X$ is smooth, and $a_{K}(T_{K})=0$ by assumption. We conclude
that $a_{\gk (x)}([\sheaf{T}\too X]_{\gk (x)})=0$ by Corollary~\ref{SpCor} above. 
\end{proof}
\end{emptythm}

\goodbreak
\section{The splitting principle}
\label{SpPrincipleSect}\bigbreak

\begin{emptythm}
\label{ReflGrSubSect}
{\bf Orthogonal reflection groups.}
We recall here and in the following three subsections some definitions and properties of
orthogonal reflection groups, merely to fix our notations. We refer to the standard
reference Bourbaki~\cite{LIE4-6} for details and more information, but see also the
book~\cite{RGInvTh} by Kane.

\smallbreak

Throughout this section we denote by~$\gk$ a field of characteristic~$\not= 2$.

\medbreak

Let~$(V,b)$ be a regular symmetric bilinear space of finite dimension over~$\gk$, and
$v\in V$ an anisotropic vector, \ie $b(v,v)\not=0$. Then
$$
s_{v}\, :\; V\,\too\, V\, ,\; w\,\longmapsto w-\frac{2b(v,w)}{b(v,v)}\cdot v\, ,
$$
is an element of the orthogonal group~$\Orth (V,b)$, called
the {\it reflection} associated with~$v$.
\end{emptythm}

\begin{emptythm}
\label{OrthReflGrDef}
{\bf Definition.}
Let~$(V,b)$ be a regular symmetric bilinear space of finite dimension over~$\gk$.
A finite subgroup of~$\Orth (V,b)$, which is generated by reflections, is called a
{\it (finite) orthogonal reflection group} over the field~$\gk$.

\medbreak

Let~$W\subset\Orth (V,b)$ be such a orthogonal reflection group.
Since~$b$ is non singular the homomorphism
$$
\hat{b}\, :\; V\,\too\, V^{\vee}\, :=\;\Hom_{\gk}(V,\gk)\, ,\; v\,\longmapsto\,
v^{\vee}\, :=\; b(\, -\, ,v)
$$
is an isomorphism. Then $b^{\vee}(v^{\vee},w^{\vee})=b(v,w)$ defines a regular
symmetric bilinear form on~$V^{\vee}$, and $v\mapsto v^{\vee}$ is an isometry
$(V,b)\xrightarrow{\simeq} (V^{\vee},b^{\vee})$. The group~$W$ acts on~$V^{\vee}$
as well via $(w.f)(x):=f(w^{-1}.x)$ for all $f\in V^{\vee}$, $w\in W~$, and~$x\in V$.

\smallbreak

Then  we have $s_{v}.f=s_{v^{\vee}}(f)$ for all $f\in V^{\vee}$ and anisotropic
vectors $v\in V$, where~$s_{v^{\vee}}$ is the reflection on~$v^{\vee}$ in~$(V^{\vee},b^{\vee})$.
Hence~$W$ is isomorphic to an orthogonal reflection group in $\Orth (V^{\vee},b^{\vee})$.

\smallbreak

We quote now two facts from Bourbaki~\cite[Chap.\ V, \S 5, Thm.\ 3 and Ex.\ 8]{LIE4-6}.
(Note for the second assertion that in a vector space with a regular bilinear form over
a field of characteristic~$\not= 2$ a pseudo reflection in the orthogonal group
is automatically a reflection, see~\cite[Chap.\ V, \S 2, no 3]{LIE4-6}.)
\end{emptythm}

\begin{emptythm}
\label{InvThm}
{\bf Theorem (Chevalley-Shephard-Todd-Bourbaki).}
{\it
Let~$W\subset\Orth (V,b)$ be a orthogonal reflection group, where~$(V,b)$ is a
regular symmetric bilinear space over the field~$\gk$. Denote by~$\SymAlg (V^{\vee})$ the
symmetric algebra of the dual space~$V^{\vee}$, and let~$f\in V^{\vee}$ be a non zero linear
form. Assume that $\khar\gk$ does not divide~$|W|$.

\smallbreak

Then:

\smallbreak

\begin{itemize}
\item[(i)]
the algebra of invariants $\SymAlg (V^{\vee})^{W}$ is a polynomial algebra; and

\smallbreak

\item[(ii)]
the isotropy group of~$f$,
$$
W_{f}\, :=\;\big\{\, w\in W\, |\, w.f=f\,\big\}
$$
is a orthogonal reflection group as well.
\end{itemize}
}

\medbreak

\noindent
Note that if~$f=v^{\vee}$ for some~$v\in V$ then we have
$$
W_{f}=W_{v}\, :=\;\big\{\, w\in W\, |\, w.v=v\,\big\}\, ,
$$
since $v^{\vee}(x)=v^{\vee}(w^{-1}.x)$ for all $x\in V$ is equivalent to
$b(w.v,x)=b(v,x)$ for all~$x\in V$ and so equivalent to $w.v=v$ since~$b$
is non singular.
\end{emptythm}

\begin{emptythm}
\label{rootSystemSubSect}
{\bf The root system of a orthogonal reflection group.}
Given such a orthogonal reflection group~$W\subset\Orth (V,b)$ let
$R_{W}$ be the set of reflections in~$W$. Recall now that
$$
w\circ s_{\alpha}\circ w^{-1}\, =\, s_{w.\alpha}\, .
$$
Hence the set~$R_{W}$ is the disjoint union of conjugacy classes
$R_{W}=\bigcup\limits_{i=1}^{m}R_{i}$. For every~$R_{i}$ we choose
an anisotropic vector $\beta_{i}$ with $s_{\beta_{i}}\in R_{i}$.
Then we have $R_{i}=\{ s_{w.\beta_{i}}|w\in W\}$ for all $1\leq i\leq m$. Let
$$
\Delta_{i}\, :=\;\big\{\, w.\beta_{i}\, |\, w\in W\,\big\}
$$
for all $1\leq i\leq m$, and set
$$
\Delta\, :=\;\bigcup\limits_{i=1}^{m}\Delta_{i}\, .
$$
Note that the sets~$\Delta_{i}$ are $W$-invariant by definition.
The set~$\Delta$ is called a {\it root system} associated with~$W$.
It has the following properties:

\smallbreak

\begin{itemize}
\item[{\bf (R1)}]
if $\alpha\in\Delta$ then $\lambda\cdot\alpha\in\Delta$ for~$\lambda\in\gk$ if and only
if~$\lambda=\pm 1$, and

\smallbreak

\item[{\bf (R2)}]
for all $\alpha,\beta\in\Delta$ we have $s_{\alpha}.\beta\in\Delta$.
\end{itemize}

\smallbreak

\noindent
(In fact, if $w.\alpha=\lambda\cdot\alpha$ then
$b(\alpha,\alpha)=b(w.\alpha,w.\alpha)=b(\lambda\cdot\alpha,\lambda\cdot\alpha)=
\lambda^{2}b(\alpha,\alpha)$, and so~$\lambda=\pm 1$, hence {\bf (R1)}. Property~{\bf (R2)}
is by construction.)

\smallbreak

\noindent
Moreover, also by construction the set $\{ s_{\alpha}|\alpha\in\Delta\}$ is the set
of all reflections in~$W$, and so in particular~$W$ is generated by
all $s_{\alpha}$, $\alpha\in\Delta$.

\smallbreak

We prove an easy lemma, which is crucial for the proof of
the main theorem.
\end{emptythm}

\begin{emptythm}
\label{VerTorLem}
{\bf Lemma.}
{\it
Let~$W\subset\Orth (V,b)$ be a orthogonal reflection group as above
and~$\Delta$ an associated root system. Let~$P_{\alpha}\in\Spec\SymAlg (V^{\vee})$
the ideal generated by~$\alpha^{\vee}$ for some~$\alpha\in\Delta$ and
$$
W_{P_{\alpha}}\, :=\;\big\{\, w\in W\, |\, w.P_{\alpha}=P_{\alpha}\,\big\}
$$
the inertia group of~$P_{\alpha}$. Then we have

\smallbreak

\begin{itemize}
\item[(i)]
$$
W_{P_{\alpha}}\, =\, W_{\pm\alpha}\ :=\;\big\{ w\in W\, |\, w.\alpha=\pm\alpha\,\big\}\,
=\, <s_{\alpha}>.W_{\alpha}\,\simeq\,\Z/2\times W_{\alpha}\, ,
$$
where~$<s_{\alpha}>=\{\id_{V},s_{\alpha}\}$ is the subgroup generated by~$s_{\alpha}$,
and

\medbreak

\item[(ii)]
if~$\khar\gk$ does not divide~$|W|$
\begin{equation}
\label{W-freeEq}
\bigcup\limits_{\alpha\in\Delta}\Ker (\alpha^{\vee})\, =\,
\bigcup\limits_{\id_{V}\not= w\in W}\Ker (w-\id_{V})\, .
\end{equation}
\end{itemize}
}

\begin{proof}
(i)~If $w.\alpha=-\alpha$ then $s_{\alpha}\circ w\in W_{\alpha}$ and so
$W_{\pm\alpha}=<s_{\alpha}>.W_{\alpha}$.

\smallbreak

For the isomorphism on the right hand side we have to show that $s_{\alpha}$
commutes with all elements of~$W_{\alpha}=\{ w\in W|w.\alpha=\alpha\}$.
This can be seen as follows. Since~$\alpha$ is anisotropic there exists a
regular subspace $H\subset V$ with $(\gk\cdot\alpha)\perp H=V$.
As $w\in\Orth (V,b)$ we have $w.h\in H$ for all $w\in W_{\alpha}$
and $h\in H$. It follows $w(s_{\alpha}(h))=w(h)=s_{\alpha}(w(h))$ for
all $h\in H$ and $w\in W_{\alpha}$. Since we have also
$w(s_{\alpha}(\alpha))=-\alpha=s_{\alpha}(w(\alpha))$ we are done.

\medbreak

(ii)~Since $\Ker\alpha^{\vee}=\Ker (s_{\alpha}-\id_{V})$ the left hand side
is contained in the right hand side.

\smallbreak

For the other direction let $x\in V\setminus\bigoplus\limits_{\alpha\in\Delta}\Ker (\alpha^{\vee})$.
By Theorem~\ref{InvThm}~(ii) above we know that $W_{x}=\{ w\in W|w.x=x\}$ is an
orthogonal reflection group. Hence if~$W_{x}\not=\{\id_{V}\}$ there exists a reflection
$s_{\alpha}\in W$, $\alpha\in\Delta$, such that $s_{\alpha}.x=x$, or equivalently $\alpha^{\vee}(x)=0$,
a contradiction.
\end{proof}

\medbreak

We state and prove now our main theorem.
\end{emptythm}

\begin{emptythm}
\label{mainThm}
{\bf Theorem.}
{\it
Let~$W$ be a orthogonal reflection group over the field~$\gk$, whose characteristic
is coprime to the order of~$W$, and~$\CM_{\ast}$ a cycle module over~$\gk$. Then
a cohomological invariant
$$
a\, :\;\HM^{1}(\, -\, ,W)\,\too\,\CM_{n}(\, -\, )
$$
over~$\gk$ is trivial if and only if its restrictions to all elementary abelian
$2$-subgroups of~$W$, which are generated by reflections, are trivial.
}

\bigbreak

For the proof we have to describe a versal $W$-torsor over~$\gk$. Let for
this~$(V,b)$ be a regular symmetric bilinear space over~$\gk$ with
$W\subset\Orth (V,b)$. Denote further by~$\Delta$ a root system
associated with~$W$.
\end{emptythm}

\begin{emptythm}
\label{ReflGrVerTorSubSect}
{\bf A versal torsor for~$W$.}
Define~$U\subset\A (V)=\Spec\SymAlg (V^{\vee})$ as in Example~\ref{VerTorExpl},
\ie~$U$ is the open complement of the union of closet sets~$\sheaf{V}_{w}$, $w\in W$,
where~$\sheaf{V}_{w}$ is the closed set defined by the ideal generated by all $f\circ (\id_{V}-w)$,
$f\in V^{\vee}$. The group~$W$ acts freely on~$U$, and the generic fiber of the quotient
morphism $q:U\too U/W$ is a versal $W$-torsor over~$\gk$, see Example~\ref{VerTorExpl}.

\smallbreak

By~(\ref{W-freeEq}) we have that
\begin{equation}
\label{W-free2Eq}
\bigcup\limits_{\alpha\in\Delta}\Ker (\id_{L}\otimes\,\alpha^{\vee})\, =\,
\bigcup\limits_{\id_{V}\not= w\in W}\Ker\big((\id_{L}\otimes\, w)-\id_{L\otimes_{\gk}V}\big)
\end{equation}
for all field extensions $L\supseteq\gk$. We get
$U=\Spec\big(\,\SymAlg (V^{\vee})[g_{\Delta}^{-1}]\,\big)$, where we have set
$$
g_{\Delta}\,  :=\;\prod\limits_{\alpha\in\Delta}\alpha^{\vee}\;\in\,\SymAlg (V^{\vee})\, .
$$
Hence the quotient morphism $q:U\too U/W$ corresponds to the embedding of rings
$$
\big(\,\SymAlg (V^{\vee})\,\big)^{W}[g_{\Delta}^{-1}]\,\too\,
\SymAlg (V^{\vee})[g_{\Delta}^{-1}]\, ,
$$
and so a versal $W$-torsor over~$\gk$, which is the generic fiber of~$q$, is
equal to the Galois extension $\Spec E\too\Spec E^{W}$ with Galois group~$W$.
Here~$E$ denotes the fraction field of~$\SymAlg (V^{\vee})$. We set in the
following $K:=E^{W}$, and denote by $[E/K]\in\HM^{1}(K,W)$ the class of the $W$-Galois
algebra $E\supset K$, which is a versal $W$-torsor over~$\gk$.
\end{emptythm}

\begin{emptythm}
\label{UnramifiedExtSubSect}
{\bf An unramified extension.}
Let~$Q$ be a prime ideal of height one in $\SymAlg (V^{\vee})^{W}$, which is not in the open
subscheme $U/W=\Spec\big(\SymAlg (V^{\vee})^{W}[g_{\Delta}^{-1}]\big)$, \ie $g_{\Delta}\in Q$.
The local ring $R:=\big(\SymAlg (V^{\vee})^{W}\big)_{Q}$ at~$Q$ is a discrete valuation ring. Let~$P$
be a prime ideal in~$\SymAlg (V^{\vee})$ above~$Q$. Then there exists~$\alpha\in\Delta$,
such that $P=P_{\alpha}=\SymAlg (V^{\vee})\cdot\alpha^{\vee}$. The Galois group~$W$ of $E\supset K$
acts transitively on the prime ideals above~$Q$, and by Lemma~\ref{VerTorLem}~(ii) we know
that the inertia group $W_{P_{\alpha}}=\{ w\in W|w.P_{\alpha}=P_{\alpha}\}$ is equal
$$
W_{\pm\alpha}\, =\,\big\{ w\in W\, |\, w.\alpha=\pm\alpha\,\big\}
\, =\, <s_{\alpha}>.W_{\alpha}\,\simeq\,\Z/2\times\; W_{\alpha}\, .
$$

\smallbreak

Denote by~$F_{\alpha}$ the fixed field of~$W_{\pm\alpha}$ in~$E$, by $\iota_{\alpha}$ the
embedding $K\subseteq F_{\alpha}$, and by~$\tilde{S}$ the integral closure of~$R$ in~$F_{\alpha}$.
We set $S_{\alpha}:=\tilde{S}_{\tilde{S}\cap P}$. This is a discrete valuation ring with maximal
ideal $Q_{\alpha}:=(\tilde{S}\cap P)\cdot\tilde{S}_{\tilde{S}\cap P}$. By construction
$S_{\alpha}\supseteq R$ is unramified, and the residue field~$\gk (Q_{\alpha})$
of~$S_{\alpha}$ is equal to the one of~$R$, which we denote~$\gk (Q)$. Hence by
the cycle module Axiom~{\bf (R3a)}, see Section~\ref{2edResMapSubSect}, we have
a commutative diagram
$$
\xymatrix{
\CM_{\ast}(F_{\alpha}) \ar[r]^-{\partial_{Q_{\alpha}}} & \CM_{\ast -1}(\gk (Q))
\\
\CM_{\ast}(K) \ar[r]^-{\partial_{Q}} \ar[u]^-{\iota_{\alpha\,\CM}} & \CM_{\ast -1}(\gk (Q))
    \ar[u]_-{=}
}
$$
for all cycle modules~$\CM_{\ast}$ over~$\gk$.
\end{emptythm}

\begin{emptythm}
\label{PfmainThmSubSect}
{\bf Proof of Theorem~\ref{mainThm}.}
The proof is by induction on $m=|W|$. If~$m=1$ or $m=2$ there is nothing
to prove, so let $m\geq 4$. Using the induction hypothesis we show first the following.

\medbreak

\noindent
{\bf Claim.}
$a_{K}([E/K])\;\in\,\CM_{n,unr}(\A (V)/W)$.

\medbreak

To prove the claim we have to show
\begin{equation}
\label{PfmainThmEq1}
\partial_{Q}\big(\, a_{K}([E/K])\,\big)\, =0
\end{equation}
for all prime ideals~$Q$ of height one in~$\SymAlg (V^{\vee})^{W}$. This is clear by
Theorem~\ref{specializationThm} if~$Q$ is in the open subset
$U/W\subset\A (V)/W=\Spec\big(\SymAlg (V^{\vee})^{W}\big)$
since~$[E/K]$ is by construction the generic fiber of $U\too U/W$.

\smallbreak

So assume~$Q\not\in U/W$. Then~$Q$ contains~$g_{\Delta}$.  Let~$Q_{\alpha}$
and~$\iota_{\alpha}:K\subseteq F_{\alpha}$ be as in Section~\ref{UnramifiedExtSubSect}
above. By the diagram at the end of Section~\ref{UnramifiedExtSubSect} is is enough to
show
$$
\partial_{Q_{\alpha}}\big(\, (\iota_{\alpha})_{\CM}(a_{K}([E/K]))\,\big)\; =\, 0\, .
$$
Since~$a$ is an invariant we have
$(\iota_{\alpha})_{\CM}(a_{K}([E/K]))=a_{F_{\alpha}}(r_{\iota_{\alpha}}([E/K]))$,
and therefore this equation is equivalent to
\begin{equation}
\label{PfmainThmEq2}
\partial_{Q_{\alpha}}\big(\, a_{F_{\alpha}}(r_{\iota_{\alpha}}([E/K]))\,\big)\; =\, 0\, .
\end{equation}

\medbreak

Now we distinguish two cases:

\smallbreak

\begin{itemize}
\item[(a)]
$F_{\alpha}=K$:
Then~$W= <s_{\alpha}>.W_{\alpha}\simeq\Z/2\times W_{\alpha}$ (in the notation of
Section~\ref{UnramifiedExtSubSect}), and therefore we have
$$
\HM^{1}(\, -\, ,W)\,\simeq\,\HM^{1}(\, -\, ,\Z/2)\times\HM^{1}(\, -\, ,W_{\alpha})\, .
$$
We claim that $a_{\ell}(x,y)=0$
for all $(x,y)\in\HM^{1}(\ell,\Z/2)\times\HM^{1}(\ell,W_{\alpha})$ and
all $\ell\in\Fields_{\gk}$. This clearly implies $a\equiv 0$, and so
$\partial_{Q}(a_{K}(T))=0$.

\smallbreak

Let for this~$\ell\in\Fields_{\gk}$, and $x\in\HM^{1}(\ell,\Z/2)$, and consider~$\Fields_{\ell}$
as full subcategory of~$\Fields_{\gk}$, \cf~Section~\ref{NotationsSubSect}. The maps
$$
b^{x}_{L}\, :\;\HM^{1}(L,W_{\alpha})\,\too\CM_{n}(L)\, ,\; z\,\longmapsto a_{L}(r_{j}(x),z)\, ,
$$
where $j:\ell\too L$ is the structure map in~$\Fields_{\ell}$, define an invariant
of degree~$n$ of~$W_{\alpha}$ over~$\ell$ with values in~$\CM_{n}$, \ie we have
$b^{x}\in\Inv^{n}_{\ell}(W_{\alpha},\CM_{\ast})$.

\smallbreak

Let~$H\subseteq W_{\alpha}$ be a $2$-subgroup generated by reflections. Then
the subgroup $H':=<s_{\alpha}>.H$ of~$W$ is a $2$-subgroup generated by reflections
as well, and therefore by assumption the restriction $\Res_{W}^{H'}(a)$ is trivial. Now
for~$L\in\Fields_{\ell}$ with structure map $\varphi_{L}:\ell\too L$ and $z\in\HM^{1}(L,H)$
we have
$$
\Res_{W_{\alpha}}^{H}(b^{x})_{L}(z)\, =\,\Res_{W}^{H'}(a)(r_{\varphi_{L}}(x),z)\, =\, 0\, .
$$
As~$z\in\HM^{1}(L,H)$ and~$L\in\Fields_{\ell}$ were arbitrary this implies
$\Res_{W_{\alpha}}^{H}(b^{x})$ is trivial. This holds for all $2$-subgroups~$H$ of~$W_{\alpha}$
generated by reflections and therefore since $W_{\alpha}$ is also a orthogonal
reflection group by Theorem~\ref{InvThm}~(ii) we have by induction~$b^{x}\equiv 0$. In
particular, we have $0=b^{x}_{\ell}(y)=a_{\ell}(x,y)=0$ as claimed. We are done in the case
$F_{\alpha}=K$.

\bigbreak

\item[(b)]
$F_{\alpha}\not= K$: Then $W\not= <s_{\alpha}>.W_{\alpha}=W_{\pm\alpha}$.
By Example~\ref{GaloisExpl} we have
\begin{equation}
\label{PfmainThmEq3}
[E/K]_{F_{\alpha}}\, =\, r_{\iota_{\alpha}}([E/K])\, =\,\theta^{1}(T')
\end{equation}
for some~$T'\in\HM^{1}(F_{\alpha},W_{\pm\alpha})$, where
$\theta:W_{\pm\alpha}\hookrightarrow W$ is the inclusion.

\smallbreak

Let~$H$ be a $2$-subgroup of~$W_{\pm\alpha}$ generated by reflections.
Then
$$
\Res_{W_{\pm\alpha}}^{H}\big(\,\Res_{W}^{W_{\pm\alpha}}(a)\,\big)\, =\,\Res_{W}^{H}(a)\, ,
$$
and since~$H$ is also a $2$-subgroup of~$W$ generated by reflections we
have $\Res_{W}^{H}(a)\equiv 0$ by our assumption. it follows that
$\Res_{W_{\pm\alpha}}^{H}\big(\,\Res_{W}^{W_{\pm\alpha}}(a)\,\big)\equiv 0$
for all $2$-subgroups of~$W_{\pm\alpha}$ which are generated by reflections.

\smallbreak

Since~$W_{\alpha}$ is a reflection group, see Theorem~\ref{InvThm}~(ii)
also $W_{\pm\alpha}$ is one, and so we can conclude by induction that
\begin{equation}
\label{PfmainThmEq4}
\Res_{W}^{W_{\pm\alpha}}(a)\;\equiv\, 0\, .
\end{equation}

\smallbreak

\noindent
We compute $\partial_{Q_{\alpha}}\big(\, a_{F_{\alpha}}(r_{\iota_{\alpha}}([E/K))\,\big)$

\smallbreak

$$
\begin{array}{r@{\; =\;}l@{\quad}l}
       & \partial_{Q_{\alpha}}\big(\, a_{F_{\alpha}}([E/K]_{F_{\alpha}})\big) & \\[3mm]
        & \partial_{Q_{\alpha}}\big(\, a_{F_{\alpha}}(\theta^{1}(T'))\,\big)
                & \mbox{by~(\ref{PfmainThmEq3})} \\[3mm]
        & \partial_{Q_{\alpha}}\big(\,\Res_{W}^{W_{\pm\alpha}}(a)_{F_{\alpha}}(T')\,\big) &
               \mbox{by definition of~$\Res_{W}^{W_{\pm\alpha}}$} \\[3mm]
        & 0 & \mbox{by~(\ref{PfmainThmEq4}).}
\end{array}
$$
Hence~(\ref{PfmainThmEq2}) holds if~$K\not= F_{\alpha}$ and we
therefore have $\partial_{Q}(a_{K}([E/K]))=0$ as claimed.
\end{itemize}

\bigbreak

\noindent
We have proven the claim, and can now finish the proof of the theorem. By the
Chevalley-Shephard-Todd-Bourbaki Theorem~\ref{InvThm}~(i)
the $\gk$-scheme $\A (V)/W=\Spec\SymAlg (V^{\vee})^{W}$ is an affine space over~$\gk$
and so by homotopy invariance of the cohomology of cycle modules, see Rost~\cite[Prop.\ 8.6]{Ro96},
we have $\CM_{n,unr}(\A (V)/W)\simeq\CM_{n}(\gk)$. Hence by the detection
Theorem~\ref{DetectionThm} the invariant~$a$ is constant. However by assumption the restriction
of~$a$ to a $2$-subgroup generated by reflections is zero, and so~$a$ has to be constant zero.

\end{emptythm}

\begin{emptythm}
\label{mainCor}
{\bf Corollary.}
{\it
Let~$W$ be as in Theorem~\ref{mainThm} a orthgonal reflection group
and~$\CM_{\ast}$ a cycle module over a field~$\gk$, whose characteristic
is coprime to~$|W|$. Let further $G_{1},\ldots ,G_{r}$ different maximal
elementary abelian $2$-subgroups generated by reflections, which represent
all such subgroups up to conjugation, \ie if~$G$ is a maximal elementary
abelian $2$-subgroup of~$W$ generated by reflections then $G=wG_{i}w^{-1}$
for some $1\leq i\leq r$ and some~$w\in W$.

\smallbreak

Then the product of restriction morphisms
$$
\big(\,\Res_{W}^{G_{i}}\,\big)_{i=1}^{r}\, :\;\Inv_{\gk}(W,\CM_{\ast})\,\too\,
\bigoplus\limits_{i=1}^{r}\Inv_{\gk}(G_{i},\CM_{\ast})^{N_{W}(G_{i})}
$$
is injective. (Recall from Example~\ref{InnAutExpl} that the image of~$\Res_{W}^{G_{i}}$
is actually in the subgroup $\Inv_{\gk}(G_{i},\CM_{\ast})^{N_{W}(G_{i})}$ for all $1\leq i\leq r$.)
}

\begin{proof}
Let~$a\in\Inv^{n}_{\gk}(W,\CM_{\ast})$ be a non trivial invariant. Then by
Theorem~\ref{mainThm} there exists an elementary abelian $2$-subgroup~$H$
of~$W$, which is generated by reflections, such that $\Res_{W}^{H}(a)\not\equiv 0$.
Let~$G$ be a maximal elementary $2$-subgroup generated by reflections,
which contains~$H$. Then there exists $1\leq i_{0}\leq r$ and~$w_{0}\in W$, such
that $w_{0}Gw_{0}^{-1}=G_{i_{0}}$. Let $H'\subseteq G_{i_{0}}$ be the image of~$H$
under the inner automorphism $\iota_{w_{0}}:g\mapsto w_{0}\cdot g\cdot w_{0}^{-1}$ of~$W$.
The morphism
$$
\iota_{w_{0}}^{\ast}\, :\;\Inv_{\gk}(H,\CM_{\ast})\,\too\,\Inv_{\gk}(H',\CM_{\ast})
$$
is an isomorphism with inverse~$\iota_{w_{0}^{-1}}^{\ast}$, see Example~\ref{InnAutExpl}
for notations, and so
$$
0\,\not\equiv\;\iota_{w_{0}}^{\ast}\big(\,\Res_{W}^{H}(a)\,\big)\, =\,
\Res_{W}^{H'}\big(\,\iota_{w_{0}}^{\ast}(a)\,\big)\, .
$$
Now $\iota_{w_{0}}^{\ast}:\Inv_{\gk}(W,\CM_{\ast})\too\Inv_{\gk}(W,\CM_{\ast})$
is the identity by~\cite[Part I, Prop.\ 13.1]{CohInv}, and therefore
$0\not\equiv\Res_{W}^{H'}(a)\, =\,\Res_{G_{i_{0}}}^{H'}\big(\,\Res_{W}^{G_{i_{0}}}(a)\,\big)$.
It follows $\Res_{W}^{G_{i_{0}}}(a)\not\equiv 0$.
\end{proof}
\end{emptythm}

\begin{emptythm}
\label{FDefRem}
{\bf Remarks.}
\begin{itemize}
\item[(i)]
For~$W$ a symmetric group the splitting principle holds more generally. The field~$\gk$
can be arbitrary, and in particular~$\khar\gk$ can divide~$|W|$. This has been
shown by Serre~\cite[Part I, Thm.\ 24.9]{CohInv}.

\smallbreak

\item[(ii)]
Let~$V$ be a finite dimensional vector space over the field~$\gk$. A $\gk$-linear
automorphism~$s$ of~$V$ is called a {\it pseudo-reflection} if~$s-\id_{V}$ has
rank~$1$. A finite subgroup~$W$ of $\Gl(V)$ is called a {\it (finite) pseudo-reflection
group} if it is generated by pseudo-reflections. One may wonder whether for such
a subgroup~$W$ of~$\Gl (V)$ the following generalization of our main
Theorem~\ref{mainThm} is true or not:

\smallbreak

\noindent
{\it
If~$\khar\gk$ and~$|W|$ are coprime then an invariant $a:\HM^{1}(\, -\, ,W)\too\CM_{n}(\, -\, )$
is constant zero if and only if its restrictions to abelian subgroups generated by pseudo-reflections
are constant zero for all cycle modules~$\CM_{\ast}$ over~$\gk$.
}

\smallbreak

\noindent
Note that under the assumption that~$\khar\gk$ does not divide~$|W|$ all pseudo-reflections in~$W$
are diagonalizable, see \eg~\cite[Prop.\ in Sect.\ 14-6]{RGInvTh}.
\end{itemize}
\end{emptythm}

\goodbreak
\section{The splitting principle of Witt- and Milnor-Witt $K$-theory invariants of orthogonal reflection groups}
\label{W-MWKInvSect}\bigbreak

\begin{emptythm}
\label{WittInvSubSect}
{\bf Witt invariants.}
We refer to the book~\cite{QHF} for details and more information about
Witt groups.

\smallbreak

Given a field~$F$  of characteristic~$\not= 2$ we denote by~$\W (F)$ the Witt group of~$F$
and by~$\FdI^{n}(F)\subset\W (F)$ the $n$th power of the fundamental ideal of~$F$ for~$n\in\Z$,
where we set $\FdI^{n}(F)=\W (F)$ if~$n\leq 0$. Fixing a base field~$\gk$ these are functors
on~$\bigFields_{\gk}$, the category of all field extensions of~$\gk$.

\smallbreak

In the following we assume $\khar\gk\not= 2$.

\smallbreak

A {\it Witt invariant} of a $\gk$-algebraic group~$G$ of degree~$n$ is a natural transformation
$$
a\, :\;\HM^{1}(\, -\, ,G)\,\too\,\FdI^{n}(\, -\,)\, ,
$$
where we consider both sides as functors on~$\bigFields_{\gk}$.

\smallbreak

We denote the set of all Witt invariants of~$G$ of degree~$n$ over~$\gk$ by
$\Inv_{\gk}^{n}(G,\FdI^{\ast})$, and set
$$
\Inv_{\gk}(G,\FdI^{\ast})\, :=\;\bigoplus\limits_{n\in\Z}\Inv_{\gk}^{n}(G,\FdI^{\ast})\, .
$$

\medbreak

For these invariants the analog of Theorem~\ref{DetectionThm}, the detection principle, holds
as well.
\end{emptythm}

\begin{emptythm}
\label{WInvDetectionThm}
{\bf Theorem.}
{\it
Let~$G$ be a linear algebraic group over the field~$\gk$, and $T\in\HM^{1}(K,G)$
a versal $G$-torsor. Then we have for and $a,b\in\Inv_{\gk} (G,\FdI^{\ast})$:
$$
a_{K}(T)\, =\, b_{K}(T)\quad\Longrightarrow\quad a\, =\, b\, .
$$
}

\medbreak

\noindent
For~$a,b$ in $\Inv_{\gk}^{0}(G,\FdI^{\ast})$ this is proven in Serre's lectures~\cite[Sect.\ 27]{CohInv},
and the same argument works also for Witt invariants of non zero degree. In fact, one can use the same
arguments as for invariants with values in Galois cohomology.

\smallbreak

Our proof for cycle modules can be copied verbatim as well. We have only to replace
the specialization map by the so called first residue map and observe that given a field~$F$
with discrete valuation~$v$ and residue field~$F$ the second residue map
$\partial_{v,\pi}:\W(F)\too\W(F(v))$ associated with some uniformizer~$\pi$ of~$v$ maps
$\FdI^{n}(F)$ into~$\FdI^{n-1}(F(v))$ for all~$n\in\Z$, see Arason~\cite[Satz 3.1]{Ar75}.

\medbreak

Using the detection principle we can now copy the proof of the splitting principle in
Section~\ref{PfmainThmSubSect} verbatim to get the following result.
\end{emptythm}

\begin{emptythm}
\label{WInvmainThm}
{\bf Theorem.}
{\it
Let~$W$ be a orthogonal reflection group over the field~$\gk$, whose characteristic
is coprime to the order of~$W$. Then a Witt invariant of degree~$n$
$$
a\, :\;\HM^{1}(\, -\, ,W)\,\too\,\FdI^{n}(\, -\, )
$$
over~$\gk$ is trivial if and only if its restrictions to all elementary abelian
$2$-subgroups of~$W$, which are generated by reflections, are trivial.
}
\end{emptythm}

\begin{emptythm}
\label{MWKInvSubSect}
{\bf Milnor-Witt $K$-theory invariants.}
The $n$th  {\it Milnor-Witt $K$-group} of a field~$F$ of characteristic not~$2$, which we denote
by~$\MWK_{n}(F)$,  can be defined using generators and relations as in Morel's book~\cite[Def.\ 3.1]{A1AlgTop},
or equivalently, see Morel~\cite{Mo04}, as the pull-back
$$
\xymatrix{
\MWK_{n}(F) \ar[r]^-{f_{n,F}} \ar[d]_-{g_{n,F}} & \MK_{n}(F) \ar[d]^-{e_{n,F}}
\\
\FdI^{n}(F) \ar[r]_-{q_{n,F}} & \FdI^{n}(F)/\FdI^{n+1}(F)\, ,
}
$$
where~$e_{n,F}$ maps the symbol $\{ a_{1},\ldots ,a_{n}\}\in\MK_{n}(F)$ to the class of the $n$-fold
Pfister form $\ll a_{1},\ldots a_{n}\gg$ and~$q_{n,F}$ is the quotient map. Considering Milnor $K$-theory,
the $n$th power of the fundamental ideal, as well as Milnor-Witt $K$-theory as functors
on the category of all field extensions of a given base field~$\gk$ then $f_{n}$ and~$g_{n}$
are natural transformations of functors.

\medbreak

Let now~$\gk$ be a field of characteristic not~$2$ and~$G$ a linear algebraic group over~$\gk$.
A {\it Milnor-Witt $K$-theory invariant} of~$G$ of degree~$n$ is a natural transformation
$$
a\, :\,\HM^{1}(\, -\, ,G)\,\too\,\MWK_{n}(\, -\, )
$$
of functors on~$\bigFields_{\gk}$. We denote the set of all Milnor-Witt $K$-theory invariants of~$G$ of
degree~$n$ by $\Inv^{n}_{\gk}(G,\MWK_{\ast})$, and set
$$
\Inv_{\gk}(G,\MWK_{\ast})\, :=\;\bigoplus\limits_{n\in\Z}\Inv^{n}_{\gk}(G,\MWK_{\ast})\, .
$$
The addition in Milnor-Witt $K$-theory induces one on~$\Inv_{\gk}^{n}(G,\MWK_{\ast})$
for all~$n\in\Z$ making the set of such invariants an abelian group.

\medbreak

Given $a\in\Inv^{n}_{\gk}(G,\MWK_{\ast})$ then $f_{n}\circ a$ and $g_{n}\circ a$ are Milnor $K$-theory-
respective Witt invariants of degree~$n$, and by the very definition of Milnor-Witt $K$-theory via
the pull-back diagram above we have that $a\equiv 0$ if and only if $f_{n}\circ a\equiv 0$ and
$g_{n}\circ a\equiv 0$. Hence we have a monomorphism
$$
\Inv^{n}_{\gk}(G,\MWK_{\ast})\,\xrightarrow{\; a\mapsto (a\circ f_{n},a\circ g_{n})\;}
\Inv^{n}_{\gk}(G,\MK_{\ast})\oplus\Inv^{n}_{\gk}(G,\FdI^{\ast})\, .
$$
Consequently we have also the splitting principle for Milnor-Witt $K$-theory
 invariants of reflection groups.
\end{emptythm}

\begin{emptythm}
\label{MWKmainThm}
{\bf Theorem.}
{\it
Let~$W$ be a orthogonal reflection group over the field~$\gk$, whose characteristic
is coprime to the order of~$W$. Then a Milnor-Witt $K$-invariant of degree~$n$
$$
a\, :\;\HM^{1}(\, -\, ,W)\,\too\,\MWK_{n}(\, -\, )
$$
over~$\gk$ is trivial if and only if its restrictions to all elementary abelian
$2$-subgroups of~$W$, which are generated by reflections, are trivial.
}
\end{emptythm}

\bibliographystyle{amsalpha}

\end{document}